\documentclass[a4paper,12pt]{article}
\usepackage[T2A]{fontenc}
\usepackage[english,russian]{babel}
\usepackage{amsfonts,amssymb,amscd,amsmath}
\usepackage[dvips]{graphicx}
\usepackage[right=25mm,left=25mm,top=30mm,
bottom=40mm]{geometry}

\begin{document}
\section*{Stability under
constantly acting perturbations
for difference equations and
averaging}
\begin{center} Vladimir Burd \\
Yaroslavl State University,
Russia\\vburd@member.ams.org 
\end{center}

1. {\bf Introduction}.
\par
In 1967 Banfi established in [3] that
uniform asymptotic stability of 
a solution of an
averaged ordinary differential
equation implies closeness of
solutions of an exact and an averaged
equations on an infinite interval 
provided the solutions
have close initial conditions. 
Similar results have been obtained 
for integro-differential equations by
Filatov in [8], for ordinary 
differential equations with 
slow and fast time by Sethna in [13],
and for functional differential 
equations of retarded type by Burd 
in [5].
\par
In this paper we consider the 
problem of closeness of
solutions of an exact and an averaged difference 
equations on an infinite interval. 
Appropriate 
assertions are derived from one 
special theorem on 
the stability under constantly 
acting perturbations. We note that
 the 
stability under constantly acting
perturbations is sometimes called 
the total stability (see, Agarwal 
[1], Hahn [9]).  

\bigskip

2. {\bf Theorem on the stability 
under constantly acting perturbations}.

\bigskip

{\bf Basic notation}.
We will use the following notation:
$|x|$ is a norm of vector 
$x\in {\cal R}^m, {\cal N}$ is the 
set of nonnegative integers,
$B_x(K)=
\{x: x\in {\cal R}^m,\,|x|\le K\},
\, G={\cal N}\times B_x(K)$.
Let $f(n,x)$ be a function that 
is defined on $G$ with values in 
${\cal R}^m$ and is bounded in the 
norm. Let $N\in{\cal N}$. Let us
assume
$$
S_{x,N}(f)=\sup_{n\in{\cal N}}
\left|\sum_{k=n}^{n+N-1}
f(k,x)\right|,\quad
x\in B_x(K).
$$

\bigskip

{\bf Lemma}. {\it Let the 
function $f(n,x)$ be continuous
in x uniformly with respect to 
$n\in{\cal N}$. 
Assume that values of the function
$x(n),(n\in{\cal N})$ belong to
$B_x(K)$.
\par
Then for any $\eta>0$ there
exists a number
$\varepsilon>0$ such that
$$
\sup_{n\in{\cal N}}
\left|\sum_{k=n}^{n+N-1}
f(k,x(k))\right|<\eta,
$$
if 
$$
S_{x,N}(f)<\varepsilon
$$}.

\bigskip

{\bf Proof}. Based on the 
conditions of the Lemma
there exists $\delta>0$ such that 
$|f(n,x_1)-f(n,x_2)|<\eta$,
if $|x_1-x_2|<\delta$.
We denote by $x^0(n)$ a 
function with values in
$B_x(K)$ that 
has a finite set of values and 
satisfies 
$|x(n)-x^0(n)|<\delta,\,
\,n\in{\cal N}$. Such 
function $x^0(n)$ exists as
$B_x(K)$ is compact. The 
function  
$x^0(n)$ has no more than 
$l$ different values, where
number $l$ depends only on
$\delta$. The statement of the
Lemma follows from the inequality
$$
\left|\sum_{k=n}^{n+N-1}f(k,x(k))
\right|
\le\left|\sum_{k=n}^{n+N-1}
[f(k,x(k))-f(k,x^0(k))]\right|
+\left|\sum_{k=n}^{n+N-1}
f(k,x^0(k))\right|.
$$

\bigskip 

3. {\bf The main theorem}. 
\par
We consider the following difference
equation in ${\cal R}^m$
$$
\Delta x(n)=X(n,x)+R(n,x),\quad
n=n_0,\,n_0+1,\dots,
\eqno(1)
$$
where $\Delta x(n)=x(n+1)-x(n)$ 
and functions
$X(n,x)$ and $R(n,x)$ 
are defined on $G$.
\par
Alongside with the equation 
(1) we consider the unperturbed
difference equation 
$$
\Delta y(n)=X(n,y). \eqno(2)
$$
\par
We suppose that equation 
(2) has a solution 
$\psi(n,n_0,\xi_0)\,
(\psi(n_0,n_0,\xi_0)=\xi_0)$ 
which is defined for all
$n\ge n_0$ and $\psi(n,n_0,\xi_0)$
together with its
$\rho$-neighborhood $(\rho>0)$ 
remains in the 
interior of the set $G$.

\bigskip

{\bf Theorem 1}. {\it Let  
function $X(n,x)$ be bounded
on $G$ and satisfies
the Lipschitz condition 
$$
|X(n,x_1)-X(n,x_2)|\le L|x_1-x_2|,
\quad x_1,x_2\in B_x(K).
\eqno(3)
$$
Let function $R(n,x)$ be
continuous in x uniformly with
respect to 
$n\in{\cal N}$ and be bounded
on $G$. Let a solution 
$\psi(n,n_0,\xi_0)$ of equation
(2) be uniformly asymptotically
stable.
\par
Then for any $\varepsilon>0\,
(0<\varepsilon<\rho)$ there 
exists  $N\in{\cal N}$ and 
numbers  $\eta_1(\varepsilon),\,
\eta_2(\varepsilon)$ such that
for all solutions 
$x(n,n_0,x_0)\in B_x(K)\,
(x(n_0,n_0,x_0)=x_0)$ of equation
(1) with initial values 
satisfying inequality 
$$
|x_0-\xi_0|<\eta_1(\varepsilon)
$$
and for all functions $R(n,x)$ 
satisfying inequality
$$
S_{x,N}(R)<\eta_2(\varepsilon)
$$
the inequality 
$$
|x(n,n_0,x_0)-\psi(n,n_0,\xi_0)|
<\varepsilon,\quad n\ge n_0
\eqno(4)
$$
holds.}

\bigskip

{\bf Proof}. Let 
$y(n,n_0,x_0)$ be a solution
of equation (2) with the same
initial condition as
the solution 
$x(n,n_0,x_0)$ of equation (1). 
These solutions satisfy the 
following equations respectively
$$
y(n,n_0,x_0)=x_0+\sum_
{k=n_0}^{n-1}X(k,y(k,n_0,x_0)),
$$
$$
x(n,n_0,x_0)=x_0+\sum_
{k=n_0}^{n-1}[X(k,x(k,n_0,x_0))
+R(k,x(k,n_0,x_0))].
$$
It follows the inequality
$$
|x(n,n_0,x_0)-
y(n,n_0,x_0)|\le\sum_{k=n_0}^{n-1}
\left|X(k,x(k,n_0,x_0))-
X(k,y(k,n_0,x_0))\right|+
$$
$$
\left|\sum_{k=n_0}^{n-1}
R(k,x(k,n_0,x_0))\right|.
$$
Using condition (3) of Theorem
we obtain  
$$
|x(n,n_0,x_0)-y(n,n_0,x_0)|
\le L \sum_{k=n_0}^{n-1}|
x(k,n_0,x_0)-y(k,n_0,x_0)|+f(n),
$$
where
$$
f(n)=\left|\sum_{k=n_0}^{n-1}
R(k,x(k,n_0,x_0))\right|=
\left|\sum_{k=n_0}^{n_0+N-1}
R(k,x(k,n_0,x_0))\right|,
$$
and $N=n-n_0$.
A well known inequality 
(see, for example, [2]) implies
$$
|x(n,n_0,x_0)-y(n,n_0,x_0)|\le 
f(n)+L\sum_{k=n_0}^{n-1}f(k)
(1+L)^{n-1-k}
$$
Therefore for $n_0\le n\le 
n_0+N$ the upper bound for  
$$
|x(n,n_0,x_0)-y(n,n_0,x_0)|
$$
depends on the values $f(n)\,(n=n_0,
\dots,n_0+N-1)$.
From uniform asymptotic stability of 
solution $\psi(n,n_0,\xi_0)$
of equation (2) follows that
there exist numbers $\delta<\varepsilon$ 
and $T\in{\cal N}$ such that
inequality
$|x_0-\xi_0|<\delta$
implies
$$
\begin{array}{l}|y(n,n_0,x_0)-
\psi(n,n_0,\xi_0)|<\frac
{\varepsilon}{2} \quad
n\ge n_0,\\|y(n_0+T,n_0,x_0)-
\psi(n_0+T,n_0,\xi_0)|<
\frac{\delta}{2}.
\end{array} \eqno(5)
$$
We now set $N=T$. The Lemma 
implies that we can find a number 
$\eta_2(\varepsilon)$
such that  
$$
|x(n,n_0,x_0)-y(n,n_0,x_0)|<
\frac{\delta}{2},\quad n_0\le n
\le n_0+T.
\eqno(6)
$$
holds. Then
$$
|x(n,n_0,x_0)-\psi_(n,n_0,
\xi_0)|<\frac{\varepsilon}{2}+
\frac{\delta}{2}<
\varepsilon,\quad n_0\le n\le 
n_0+T.
$$
Furthermore, (5) and (6) imply
$$
|x(n_0+T,n_0,x_0)-\psi(n_0+T,n_0,
\xi_0)|<\delta.
$$
Hence for the interval 
$[n_0,n_0+T]$ the solution 
$x(n,n_0,x_0)$
remains in $\varepsilon$-neighborhood 
of the solution 
$\psi(n,n_0,\xi_0)$ and
at the moment $n=n_0+T$ 
belongs to the $\delta$-neighborhood
of $\psi(n,n_0,\xi_0)$. 
\par
We now consider $n=n_0+T$ as an 
initial moment.  Using the same 
arguments as above we obtain
$$ 
|x(n,n_0,x_0)-\psi(n,n_0,\xi_0)|
\le\varepsilon, \quad n_0+T\le
n\le n_0+2T
$$
and
$$
|x(n_0+2T,n_0,x_0)-\psi(n_0+2T,
n_0,\xi_0|\le\delta.
$$ 
Repetitive application 
of the same argument
completes the proof of the 
Theorem.

\bigskip

We note that last part of proof
of the Theorem 1 uses the
reasoning similar to the Lemma 6.3
from [4].

\bigskip

{\bf Remark 1}. The statement of
Theorem 1 differs from the statements
of known theorems on the stability 
under
constantly acting perturbations 
[1, 9 - 12] 
in using a more general assumption
on the 
perturbation $R(n,x)$. Theorem 1 
implies Halanay's theorem, if number
$\eta_2(\varepsilon)=N\delta_2(
\varepsilon)$, where 
$\delta_2(\varepsilon)$
is a number from definition 5.13.1 
[1]. If we assume 
$$
S_{x,N}(f)=\sup_{n\in{\cal N}}
\sum_{k=n}^{n+N-1}
|f(k,x)|,\quad
x\in B_x(K),
$$
then we obtain a difference
analog of the Theorem 24.1 from [11].

\bigskip

{\bf Remark 2}. We start with the 
following definition.

\bigskip
 
{\bf Definition}. The solution 
$\psi(n,n_0,\xi_0)$ is called
uniformly asymptotically stable 
with respect to a part of the
variables 
$\psi_1,\dots,\psi_k,\,k<m$,
if its asymptotically stable
in the sense of Lyapunov with 
respect to a part of the variables 
$\psi_1,\dots,\psi_k,\,k<m$ and 
if for any number 
$\gamma>0$ there exists 
a number
$T(\gamma)\in{\cal N}$,
such that for the solution 
$y(n,n_0,x_0)$ of the equation
(2) is satisfied inequality
$$
|y_i(n,n_0,x_0)-\psi_i(n,n_0,
\xi_0)|<\gamma,\quad  n\ge
n_0 + T(\gamma),
\quad i=1,\dots,k,
$$
for any initial moment 
$n_0$ and any initial values
$x_0$ from the domain of the 
attraction of solution 
$\psi(n,n_0,\xi_0)$ with 
respect to a part of the variables
(i.e. from domain where is 
satisfied the limit equality 
$$
\lim_{n\to\infty}|y_i(n,n_0,x_0)-
\psi_i(n,n_0,\xi_0)|=0,\quad 
i=1,\dots,k.)
$$

\bigskip

If $k=m$, the definition above 
coincides with the definition of 
uniform asymptotic 
stability (see, for example, [6]).
A detailed discussion of the 
stability theory
with respect to a part of the variables
is given, for example, in [14].
\par
If the solution
$\psi(n,n_0,\xi_0)$ of equation
(2) uniformly asymptotically
stable only with respect to 
a part of the variables 
$\psi_1,\dots,\psi_k,\, k<m$, 
then inequality (4) 
in the statement of the Theorem 1
can be replaced with the inequality
$$
|x_i(n,n_0,x_0)-\psi_i(n,n_0,
\xi_0)|<\varepsilon, \quad i=1,
\dots,k.
$$

\bigskip

4. {\bf Averaging on an 
infinite interval}.
\par
Theorem 1 is applicable to the 
problem of averaging on 
an infinite interval for 
difference equations.
\par
We consider the following difference
equation in ${\cal R}^m$
$$
\Delta x(n)=\varepsilon 
X(n,x), \eqno(7)
$$
where $\varepsilon>0$ is a small
parameter, $X(n,x)$ is
defined for $(n,x)\in G$.

\bigskip

{\bf Theorem 2}. {\it Let 
\newline
1) function $X(n,x)$ be 
continuous in $x$
uniformly with respect to 
$n\in{\cal N}$;\newline
2) $\left|X(n,x)
\right|\le M_1<\infty, 
\quad (n,x)
\in G$;\newline
3) the limit 
$$
\lim_{N\to\infty}\frac{1}{N}
\sum_{k=n}^{n+N-1} X(k,x)
=\bar{X}(x)
$$
exists uniformly with 
respect to $n$ for any $(n,x)
\in G$ and 
$\bar{X}(x)$ be bounded
in the norm
$$
|\bar{X}(x)|\le M_2<\infty,
\quad x\in B_x(K);
$$
\newline
4) function $\bar{X}(x)$
satisfies the Lipschitz 
condition
$$
\left|\bar{X}(x_1) -
\bar{X}(x_2)\right|
\le L|x_1-x_2|,
\quad x_1,x_2\in B_x(K),\,
$$
\newline
5) the averaged equation
$$
\Delta x(n)=\varepsilon
\bar{X}(x) 
$$
has a uniformly asymptotically
stable solution $\psi(n,n_0,
\xi_0)$
(uniformly asymptotically
stable with respect to 
a part of the variables
$x_1,\dots,x_k,\,k<m$),
which with its 
$\rho$-neighborhood ($\rho>0$) 
belong to $G$.
\par
Then for any $\alpha\,
(0<\alpha<\rho)$ there exists
$\varepsilon_1(\alpha)\,
(0<\varepsilon_1<\varepsilon_0)$
and  $\beta(\alpha)$ such that
for all
$0<\varepsilon<\varepsilon_1$ 
the solution 
$\varphi(n,n_0,x_0)\in B_x(K)$ of 
equation (7) with an initial
condition satisfying inequality
$$ 
|x_0-\xi_0|
<\beta(\alpha)\quad
(|x_{0i}-\xi_{0i}|<
\beta(\alpha),\, i=1,\dots,k<m)
$$
we have
$$
|\psi(n,n_0,\xi_0)-
\varphi(n,n_0,x_0)|<\alpha,
\quad n\ge n_0
$$
$$
(|\psi_i(n,n_0,\xi_0)-
\varphi_i(n,n_0,x_0)|<
\alpha,\quad i=1,\dots,k<m\quad
n\ge n_0).
$$
}

\bigskip

The equation (7) can be written
in the form
$$
\Delta x(n)=\varepsilon
\bar{X}(x)+\varepsilon
R(n,x),
$$
where $R(n,x)=X(n,x)-\bar{X}(x)$. 
We show that Theorem 2 follows
from Theorem 1. Given $\alpha$ 
we choose
$N(\alpha)=[\frac{1}
{\varepsilon}]$ where $[x]$ is the
integer part of $x$.
Then 
$$
\left|\varepsilon\sum_{k=n}^{n+
[\frac{1}{\varepsilon}]-1}R(k,x)
\right|\le\frac{1}{[\frac{1}
{\varepsilon}-1]}\left|
\sum_{k=n}^{n+
[\frac{1}{\varepsilon-1}]}R(k,x)
\right|.
$$ 
Therefore condition 3)
of the Theorem 2 implies that
for a sufficiently small
$\varepsilon$
the function $R(n,x)$  
satisfies of conditions of 
Theorem 1.

\bigskip

4. {\bf Averaging on an infinite
interval of systems with the
right-hand side that vanishes over 
time}.
\par
We now consider the following
difference equation in ${\cal R}^m$
$$
\Delta x(n)=\frac{1}{n} 
X(n,x),\quad n=n_0,n_0+1,\dots,
\eqno(8)
$$
where  $X(n,x)$ is
defined for $(n,x)\in G$.

\bigskip

{\bf Theorem 3}. {\it Let 
\newline
1) function $X(n,x)$ be 
continuous in $x$
uniformly with respect to 
$n\in{\cal N}$;\newline
2) $\left|X(n,x)
\right|\le M_1<\infty, 
\quad (n,x)
\in G$;\newline
3) the limit 
$$
\lim_{N\to\infty}\frac{1}{N}
\sum_{k=n}^{n+N-1} X(k,x)
=\bar{X}(x)
$$
exist uniformly with 
respect to $n$ for any $(n,x)
\in G$ and 
$\bar{X}(x)$ be bounded
in the norm
$$
|\bar{X}(x)|\le M_2<\infty,
\quad x\in B_x(K);
$$
\newline
4) function $\bar{X}(x)$
satisfies the Lipschitz 
condition
$$
\left|\bar{X}(x_1) -
\bar{X}(x_2)\right|
\le L|x_1-x_2|,
\quad x_1,x_2\in B_x(K),\,
$$
\newline
5) the averaged equation
$$
\Delta x(n)=\frac{1}{n}
\bar{X}(x) 
$$
has a uniformly asymptotically
stable solution $\psi(n,n_0,
\xi_0)$
(uniformly asymptotically
stable with respect to 
a part of variables
$x_1,\dots,x_k,\,k<m$),
which with its 
$\rho$-neighborhood 
($\rho>0$) is contained
in $G$.
\par
Then for any $\alpha\,
(0<\alpha<\rho)$ there exists
$\varepsilon_1(\alpha)\,
(0<\varepsilon_1<\varepsilon_0)$
and  $\beta(\alpha)$ such that
for all
$0<\varepsilon<\varepsilon_1$ 
the solution 
$\varphi(n,n_0,x_0)\in B_x(K)$ of 
equation (8), with initial
condition satisfying inequality
$$ 
|x_0-\xi_0|
<\beta(\alpha)\quad
(|x_{0i}-\xi_{0i}|<
\beta(\alpha),\, i=1,\dots,k)
$$
we have
$$
|\psi(n,n_0,\xi_0)-
\varphi(n,n_0,x_0)|<\alpha,
\quad n\ge n_0
$$
$$
(|\psi_i(n,n_0,\xi_0)-
\varphi_i(n,n_0,x_0)|<
\alpha,\quad i=1,\dots,k<m\quad
n\ge n_0).
$$
}

\bigskip

The proof of the Theorem 3 is 
quite similar to the proof of
the Theorem 2.

\bigskip

5. {\bf Dynamics of selection
of genetic population in a 
varying environment}.
\par
As an example we consider dynamics
of a selection of a Mendelian 
population with a genetic pool
made of only two alleles, that
we'll call $A$ and $a$ . We assume 
that the fitness of the genotypes $AA$, 
$Aa$, $aa$ are 
$1-\varepsilon\alpha(n),\,1,\,
1-\varepsilon\beta(n)$ respectively.
Here n is number of the generation,
$\varepsilon>0$ is a small 
parameter, $\alpha(n),\,\beta(n)$
are periodic functions 
with period $l\in{\cal N}$ and 
positive mean values. 
Let $p_n,\,q_n$ be the frequencies
of alleles $A$, $a$ in generation $n$
respectively. The evolution equation 
has the form (see, for example, [7])
$$
\Delta p_n=\varepsilon p_n(1-p_n)
\frac{\beta(n)-(\alpha(n)+\beta(n))
p_n}{1-\varepsilon[(\alpha(n)+
\beta(n))p_n^2+2\beta(n)p_n-
\beta(n)]}. \eqno(9)
$$
The averaged equation
$$
\Delta\bar{p}_n=\varepsilon
\bar{p}_n(1-\bar{p}_n)
(\beta_0-(\alpha_0+\beta_0)
\bar{p}_n)
$$
has a unique asymptotically stable
equilibrium 
$$
\bar{p}=\frac{\beta_0}
{\alpha_0+\beta_0},
$$
where $\alpha_0,\,\beta_0$ are
mean value of periodic
functions $\alpha(n),\,\beta(n)$
accordingly.
\par
Then for a sufficiently small $\varepsilon$
equation (9) has an asymptotically
stable periodic solution with
period $l\in{\cal N}$ (see [10]). 
Theorem 2 implies that for any $\delta>0$ 
there exists $\eta(\delta)$ such 
that the solution 
$p_n(0,x_0)$ of 
equation (9), with  
initial condition satisfying 
inequality
$$ 
|x_0-\xi_0|
<\eta(\delta),
$$
where $\xi_0>0,\,\xi_0\ne\frac
{\beta_0}{\alpha_0+\beta_0}$
the inequality
$$
|p_n(0,x_0)-\bar{p}_n(0,
\xi_0|<\delta,
\quad n\ge 0
$$
holds.
      
\bigskip

{\bf References}

\begin{enumerate}
\item
{R. Agarwal}, Difference equations 
and inequalities: theory, methods,
and applications, Marcel Dekker,
New York, 2000. 

\item
{R. Agarwal, M. Bohner, A.
Peterson}, Inequalities on time 
scales: A survey, Mathematical
inequalities and Applications,
4(4), (2001), 535--557.
\item
{C. Banfi}, 
Sull'approssimazione di 
processi non stazionari in 
mecanica non lineare, Bolletino
della Unione Matematica 
Italiana, 22 (1967), 442--450.
\item
{E.A. Barbashin}, Vvedenie v teoriyu
ustoichivosti, Nauka,
Moscow, 1967 (in Russian).
\item
{V.Sh. Burd}, Stability under
constantly acting disturbances
and principle of averaging
on an infinite interval for 
systems with time lag, 
Functional Differential 
equations, 4 (1997), n. 3-4,
257--264.
\item
{S. Elaydi}, An introduction to
difference equations. 
Springer-Verlag, New York, 1996.
\item
{W.J. Ewens}, Mathematical 
population genetics, 
Springer-Verlag, Berlin; New York,
1979.
\item
{A.N. Filatov}, Methods of
averaging in differential and
integro-differential equations,
Fan, Tashkent, 1971.
\item 
{W. Hahn}, Theory and applications
on Liapunov's direct method, 
Prentice-Hall, Englwood Cliffs,
N.J., 1963. 
\item
{A. Halanay, D. Wexler}, Teoria
calitativa a sistemelor cu
impulsuri, Editura Academiei
Republicii Socialiste Romania,
Bucuresti, 1968.
\item
{N.N. Krasovskii}, Stability of 
motion, Stanford University Press,
Stanford, California, 1963.
\item
{I.G. Malkin}, Theory of stability
of motion, Atomic Energy Commission,
Translation No. 3352, Washington,
D.C., 1958.
\item
{P.R. Sethna}, Systems with 
fast and slow time, The 5th
international conference on
nonlinear oscillations, 
Ukrainian Academy of Sciences,
Kiev, 1 (1971), 505--521.
\item
{V.I. Vorotnikov}, Partial stability
and control, Birkh\"auser, Boston, 
1997.
\end{enumerate}
\end{document}